\documentclass[11pt,a4paper]{amsart}
\usepackage{mathrsfs}
\usepackage{syntonly}
\usepackage{amsmath}
\usepackage{amsthm}
\usepackage{amsfonts}
\usepackage{amssymb}
\usepackage{latexsym}
\usepackage{amscd,amssymb,amsopn,amsmath,amsthm,graphics,amsfonts,mathrsfs,accents,enumerate,verbatim,calc}
\usepackage[dvips]{graphicx}
\usepackage[colorlinks=true,linkcolor=red,citecolor=blue]{hyperref}
\usepackage[all]{xy}
\usepackage{enumitem}

\date{}
\allowdisplaybreaks[4] \footskip=15pt
\renewcommand{\uppercasenonmath}[1]{}

\numberwithin{equation}{section} \theoremstyle{plain}
\newtheorem{lem}{Lemma}[section]
\newtheorem{cor}[lem]{Corollary}
\newtheorem{prop}[lem]{Proposition}
\newtheorem{thm}[lem]{Theorem}

\newtheorem{Ex}[lem]{Example}
\newtheorem{Quest}[lem]{Question}
\newtheorem{Property}[lem]{Property}
\newtheorem{Properties}[lem]{Properties}
\newtheorem{Subprops}{}[lem]
\newtheorem{Para}[lem]{}

\newtheorem{rem}[lem]{Remark}

\newtheorem*{ack*}{ACKNOWLEDGEMENTS}

\newcommand{\pf}{\noindent\begin {proof}}
\newcommand{\epf}{\end{proof}}

\newcommand{\mto}{\rightarrowtail}
\newcommand{\eto}{\twoheadrightarrow}

\newcommand{\s}{\stackrel}

\newcommand{\ME}{{\rm ME}}

\newcommand{\Ext}{{\rm Ext}}

\newcommand{\Hom}{{\rm Hom}}

\newcommand{\mcA}{\mathcal{A}}

\newcommand{\mcE}{\mathcal{E}}

\newcommand{\mcI}{\mathcal{I}}
\newcommand{\mcJ}{\mathcal{J}}

\newcommand{\mcM}{\mathcal{M}}

\newcommand{\mcT}{\mathcal{T}}
\newcommand{\mcX}{\mathcal{X}}
\newcommand{\mcD}{\mathcal{D}}
\newcommand{\mcY}{\mathcal{Y}}

\newcommand{\Arr}{{\rm Arr}}
\begin{document}
\begin{center}
{ \bf \MakeUppercase{Ideal approximation theory in Frobenius categories}}

\vspace{0.5cm}  Dandan Sun, Zhongsheng Tan, Qikai Wang and Haiyan Zhu\footnote{Corresponding author.\\
\indent Haiyan Zhu was supported by  the NSF of China (12271481).}
\end{center}
\bigskip
\centerline { \bf  Abstract}
\leftskip10truemm \rightskip10truemm \noindent
Let $\mcA$ be a Frobenius category and $\omega$ the full subcategory consisting of projective objects. The relations between special precovering (resp., precovering) ideals in $\mcA$ and special precovering (resp., preenveloping) ideals in the stable category $\mcA/\omega$ are explored. In combination with a result due to Breaz and Modoi, we conclude that every precovering or preenveloping ideal $\mcI$ in $\mcA$ with $1_{X}\in{\mcI}$ for any $X\in{\omega}$ is special. As a consequence, it is proved that an ideal cotorsion pair $(\mcI,\mcJ)$ in $\mcA$ is complete if and only if $\mcI$ is precovering if and only if $\mcJ$ is preenveloping. This leads to an ideal version of the Bongartz-Eklof-Trlifaj Lemma in $\mcA/\omega$, which states that an ideal cotorsion pair in $\mcA/\omega$ generated by a set of morphisms is complete. As another consequence, we provide some partial answers to the question about the completeness of cotorsion pairs posed by Fu, Guil Asensio, Herzog and Torrecillas.
\leftskip10truemm \rightskip10truemm \noindent
\\[2mm]

\leftskip0truemm \rightskip0truemm

\section{Introduction}

The approximation theory is one of the efficient tools to study the complicated objects by some simple objects in a category. It is an important tool in study of various categories. For example, they are used to construct resolutions and to do homological algebra: in module theory the existence of injective envelopes, projective precovers and flat covers is often used for defining derived functors, for dealing with invariants as (weak) global dimension etc.
The central role in approximation theory for the case of modules, or more general abelian or exact categories, is played by the notion of complete cotorsion pairs \cite{A6}.  Recall that the concept of cotorsion pairs was introduced by Salce \cite{A10} in the setting of the category of abelian groups. Cotorsion pairs have been used to study approximation theory in sheaf categories \cite{A4}, general Grothendieck categories \cite{A7,A8}, and more abstract exact categories \cite{A11}.

To replace objects and subclasses by morphisms and ideals, Fu, Asensio, Herzog and Torrecillas introduced ideal cotorsion pairs and the ideal approximation theory \cite{FGHT}.  In this theory, the role of the objects and subcategories in classical approximation theory is replaced by the morphisms and ideals of the category. In the meantime, they established the connection between special precovering ideals and complete ideal cotorsion pairs in exact categories. As a generalization of the classical approximation theory, the ideal approximation theory has been extended by Breaz and Modoi \cite{Breaz} to triangulated categories.

An interesting observation in \cite{Breaz} is that every precovering or preenveloping ideal in a triangulated category is special. However, in general, it is not true if we replace the triangulated category with an exact category. For example,
let $R=k[x_{1},x_{2}]$ be a polynomial ring over a field $k$ and $\mathcal{I}$ an ideal of morphisms that factor through an injective $R$-module, whence $\mathcal{I}$ is a precovering ideal in the category of $R$-modules, but it is not special since $R$ is not an injective $R$-module and does not have an epimorphism $E\rightarrow R$ with $E$ an injective $R$-module. The following is our main result which provides sufficient conditions for every precovering or preenveloping ideal to be special in an exact category.

\begin{thm}\label{ThA} Assume that $(\mcA; \mcE)$ is a Frobenius category and $\omega$ is the full subcategory consisting of projective objects. Denote by $\mcA/\omega$ the stable category of $(\mcA; \mcE)$.
\begin{enumerate}
\item If $\mcI$ is an ideal in $(\mcA; \mcE)$ such that $\omega$ $\subseteq$ {\rm Ob}($\mathcal I$), then the following statements are equivalent.
\begin{enumerate}
\item  $\mcI$ is special precovering in $\mcA$.
\item  $\mcI$ is precovering in $\mcA$.
\item  $(\mcI,\mcI^{\perp})$ is a complete cotorsion pair in $\mcA$.
\item  $(\mcI/\omega,\mcI^{\perp}/\omega)$ is a complete cotorsion pair in $\mcA/\omega$.
\end{enumerate}

\item If $\mcJ$ is an ideal in $(\mcA; \mcE)$ such that $\omega$ $\subseteq$ {\rm Ob}($\mathcal J$), then the following statements are equivalent.
\begin{enumerate}
\item  $\mcJ$ is special preenveloping in $\mcA$.
\item  $\mcJ$ is preenveloping in $\mcA$.
\item  $(^{\perp}\mcJ, \mcJ)$ is a complete cotorsion pair in $\mcA$.
\item  $(^{\perp}\mcJ/\omega, \mcJ/\omega)$ is a complete cotorsion pair in $\mcA/\omega$.
\end{enumerate}
\end{enumerate}
\end{thm}

We note that the precise definitions of (special) precovering ideals, (special) preenveloping ideals, and complete ideal cotorsion pairs can be seen in Section 2.3 below. The strategy of our proof of Themorem \ref{ThA} is to explore the relations between special precovering (resp., preenveloping) ideals in $(\mcA; \mcE)$ and special precovering (resp.,  preenveloping) ideals in $\mcA/\omega$ (see Proposition \ref{thm1}).

An immediate consequence of Theorem \ref{ThA} yields the following corollary.
\begin{cor}\label{cor3} Let $(\mcI,\mcJ)$ be an ideal cotorsion pair in a Frobenius category $(\mcA; \mcE)$. Then the following statements are equivalent.
\begin{enumerate}
\item  $\mcI$ is precovering in $\mcA$.
\item  $\mcJ$ is preenveloping in $\mcA$.
\item  $(\mcI,\mcJ)$ is a complete cotorsion pair in $\mcA$.
\item  $(\mcI/\omega, \mcJ/\omega)$ is a complete cotorsion pair in $\mcA/\omega$.
\end{enumerate}
\end{cor}

As an application of Corollary \ref{cor3}, we prove an ideal version of the Bongartz-Eklof-Trlifaj Lemma \cite{FHHZ} in the stable category $\mcA/\omega$ of a Frobenius category $(\mcA; \mcE)$ with exact coproducts, which implies that an ideal cotorsion pair in $\mcA/\omega$ generated by a set of morphisms is complete (see Corollary \ref{cor:BKT-lemma}). Examples of stable categories of Frobenius categories include but are not limited to: the derived category of a ring $R$, the homotopy category of complexes of modules and so on. We refer to \cite{Keller2006} and \cite{Krause2007} for a more detailed discussion on this matter.

For a class $\mathcal{X}$ and an ideal $\mcI$ of $(\mcA; \mcE)$, denote by $\mcI(\mathcal{X})$ the ideal of those morphisms that factor through an object in $\mathcal{X}$ and by ${\rm Ob}(\mcI):=\{X\in{\mcA}~|~1_{X}\in{\mcI}\}$. If $\mcI=\mcI({\rm Ob}(\mcI))$, then we call $\mcI$ is an object ideal. Moreover, let $\mcX$ and $\mcY$ be additive subcategories of $(\mcA; \mcE)$, the classical situation of a complete cotorsion pair $(\mcX, \mcY)$ in $(\mcA; \mcE)$ yields a complete ideal cotorsion pair $(\mcI(\mcX), \mcI(\mcY))$. However, the converse problem, due to Fu, Asensio, Herzog and Torrecillas \cite[Question 29]{FH}, that is if the cotorsion pair $({\rm Ob}(\mcI),{\rm Ob}(\mcJ))$ is necessarily complete provided that $\mcI$ and $\mcJ$ are object ideals in $\mcA$ such that $(\mcI,\mcJ)$ is a complete ideal cotorsion pair, is still open.

Let $\mcX$ be a full subcategory of an additive category $\mcD$.  We set \begin{center}{${\rm smd}(\mcX):=\{M~|~M~\mathrm{is}~\mathrm{isomorphic}~\mathrm{to}~\mathrm{a}~\mathrm{direct}~\mathrm{ summand}~ \mathrm{of}~ \mathrm{an}~\mathrm{object}~\mathrm{of}~ \mcX~\mathrm{in}~ \mcD\}.$}
\end{center}
Recall that an additive category is called a \emph{Krull-Schmidt category} if every object decomposes into a finite direct sum of objects having local endomorphism rings.  One can refer to a standard reference such as \cite{Krause} for further details on Krull-Schmidt categories.

Another consequence of Theorem \ref{ThA} yields the following result, which gives some partial answers to the Question 29 in \cite{FH}.

\begin{cor}\label{corollary2} Let $(\mcI,\mcJ)$ be a complete ideal cotorsion pair in a Frobenius category $\mcA$ such that both of the ideals $\mcI$ and $\mcJ$ are object ideals. Set $\mcX:={\rm Ob}(\mcI)$ and $\mcY:={\rm Ob}(\mcJ)$. Then the following statements are equivalent.
\begin{enumerate}
\item $(\mcX,\mcY)$ is a complete cotorsion pair in $\mcA$.
\item For every object $\underline{M}$ in ${\rm smd}(\underline{\mcX}\oplus\underline{\mcY})$, there is a triangle
$\xymatrix{\underline{M}\ar[r]&\underline{Y}\ar[r]&\underline{X}\ar[r]&\underline{M}[1]}$  in $\mcA/\omega$ with $Y\in{\mcY}$ and $X\in{\mcX}$.
\item For every object ${M}$ in ${\rm smd}({\mcX}\oplus{\mcY})$, there is a conflation
$\xymatrix{{M}\ar[r]&{Y}\ar[r]&{X}}$  in $\mcA$ with $Y\in{\mcY}$ and $X\in{\mcX}$.
\end{enumerate}

If in addition $\mcA/\omega$ is Krull-Schmidt, then any of the above hold.
\end{cor}




The contents of this paper are outlined as follows. In Section 2, we fix notations and recall some definitions and basic facts used throughout the paper. In Section 3, we first study relations between special precovering (resp., preenveloping) ideals in a Frobenius category $(\mcA; \mcE)$ and special precovering (resp.,  preenveloping) ideals in its stable category $\mcA/\omega$ (see Proposition \ref{thm1}), and then we give the proof of Theorem \ref{ThA}. Some applicaitons of Theorem \ref{ThA} are given in Section 4, including the proof of Corollary \ref{corollary2}.

When finishing the work on our paper, we learned about a bijective correspondence between complete ideal cotorsion pairs in a Frobenius category and complete ideal torsion pairs in this stable category has been obtained by Zhu, Fu, Herzog and Schlegel \cite{ZFHS}. In a forthcoming paper \cite{STWZ}, we will study ideal approximation theory in Frobenius extriangulated categories.

\section{Preliminary}

The assumptions, the notation, and the definitions from this section will be used throughout the
paper.  For more details, the reader can consult \cite{Breaz,FGHT,FH,A6}.

\subsection{Exact categories}

Recall from \cite{TBh10,Keller,Quillen} that an \emph{exact category} is a pair $(\mcA; \mcE),$  where $\mcA$ is an additive category and $\mcE$ is a class of ``short exact sequences": That is,
an actual kernel-cokernel pair
$$\Xi \colon \xymatrix@C=30pt@R=30pt{A~\ar@{>->}[r] & B \ar@{->>}[r] & C.}$$ In what follows, we call such a sequence a \emph{conflation}, and call $A\mto B$ (resp., $B\eto C$) an \emph{inflation} (resp., \emph{deflation}). Many authors use the alternate terms \emph{admissible exact sequence}, \emph{admissible monomorphism} and \emph{admissible epimorphism}.
The class $\mathcal{E}$ of conflations must satisfy exact axioms, for details, we refer the reader to \cite[Definition 2.1]{TBh10}, which are inspired by the properties of short exact sequences in any abelian category. We often write $\mcA$ instead of $(\mcA; \mcE)$ when we consider only one exact structure on $\mathcal{A}$.

Recall that an object $P$ in $\mcA$ is called \emph{projective} provided that any admissible epimorphism ending at $P$ splits.
The exact category $\mcA$ is said to have \emph{enough projective objects} provided that each object $X$ fits into a deflation $d: P \rightarrow X$ with $P$ projective. Dually, one has the notions of injective
objects and enough injective objects.

An exact category $(\mcA; \mcE)$ is said to be \emph{Frobenius} provided that it has enough projective and enough injective objects, and that the class of projective objects coincides with the class of injective objects. Denote by $\omega$ the full subcategory of $\mcA$ consisting of projective objects.
The importance of Frobenius categories lies in that they give rise
naturally to triangulated categories; see \cite{Happel}.



\subsection{Precovering and preenveloping ideals in additive categories}
An {\em ideal} $\mcI$ of an additive category $\mcA$ is an additive subbifunctor of $\Hom: \mcA^{{\rm op}} \times \mcA \to {\rm Ab};$ it associates to a pair $(A,B)$ of objects in $\mcA$ a subgroup $\mcI(A,B) \subseteq \Hom (A,B)$
and satisfies the usual conditions in the definition of an ideal in ring theory: for any morphism $g: A \to B$ in $\mcI,$ the composition $fgh$ belongs to $\mcI (X,Y),$ for any morphisms $f: B \to Y$ and $h: X \to A$ in $\mcA.$

If $A\in{\mcA}$ is an object, we say that $A$ belongs to an ideal $\mcI$ if the identity morphism $1_{A}$ belongs to
$\mcI(A,A),$ and denote by $\textrm{Ob}(\mcI)\subseteq \mcA$ the full subcategory of objects of $\mcI$.

Let $\mcI$ be a class of morphisms in an additive category $\mcA$. An $\mcI$-precover of an object $X\in{\mcA}$ is a morphism $f:A\rightarrow X$ in $\mcI$ such that any other morphism $f':A'\rightarrow X$ in $\mcI$ factors through $f$.
 The class $\mcI$ is said to be \emph{precovering} if every object in $\mcA$ has an $\mcI$-precover. The notion of $\mcI$-preenvelope is defined dually.

\subsection{Special precovering and special preenveloping ideals in exact categories} Two morphisms $a$ and $b$ in an exact category $\mcA$ are \emph{$\Ext$-orthogonal} if the morphism $\Ext(a,b) : \Ext(A_1, B_0) \to \Ext(A_0, B_1)$ of abelian groups is $0.$
If $\mcM$ is a class of morphisms in $\mcA,$ then $${^{\perp}}\mcM=\{a\ |\ \Ext(a,m)=0\ \text{for any}\ m\in\mcM\}\ \text{and}\ \mcM^{\perp}=\{b\ |\ \Ext(m,b)=0\ \text{for any}\ m\in\mcM\}$$ are ideals. One defines
an \emph{ideal cotorsion pair} consists of ideals $(\mcI, \mcJ)$ for which $\mcI = {^{\perp}}\mcJ$ and $\mcJ^{\perp} = \mcI.$

If $\mcA$ is endowed with an exact structure $\mcE$, then an $\mcI$-precover $f:A\rightarrow X$ of $X$ is said to be \emph{special} if $f$ is a deflation in $\mcE$ and it fits in a pushout diagram
$$\xymatrix@C=30pt@R=30pt{
K ~\ar@{>->}[r] \ar[d]^{g} & A' \ar@{->>}[r] \ar[d] & X \ar@{=}[d] \\
E~\ar@{>->}[r] & A \ar@{->>}[r]^{f} & X}
$$
with $\Ext(i,g)=0$ for every $i\in{\mcI}$. The class $\mcI$ is said to be \emph{special precovering} if every object in $\mcA$ has a special $\mcI$-precover. The notion of special $\mcI$-preenvelope is defined dually.

Moreover, an ideal cotorsion pair $(\mcI,\mcJ)$ is \emph{complete} if every object in $\mcA$ has a special $\mcJ$-preenvelope and a special $\mcI$-precover (see \cite{FGHT}).

\subsection{Special precovering and special preenveloping ideals in triangulated categories}
Let $\mcT$ be a triangulated catgeory. Recall from \cite{Breaz} that a morphism $f : X \to A$ in $\mcT$ is called \emph{left orthogonal} to a morphism $g : B \to Y$ in $\mcT$ and we denote this by $f \perp g$, if $g[1]{\phi}f = 0$ for all morphisms $\phi:A\to B[1]$. If $\mcM$ is a class of morphisms in $\mcT,$ then ${^{\perp}}\mcM$ is the ideal of morphisms left $\Ext$-orthogonal to every $m \in \mcM.$ One defines ${\mcM}^{\perp}$ dually.

Let $\mcI$ be an ideal in $\mcT$. Thanks to \cite[Remark 3.3.1]{Breaz}, a morphism $i : X \to A$ is a \emph{special $\mcI$-precover} if there exists a homotopy pushout diagram
$$\xymatrix@C=30pt@R=30pt{Y \ar[r] \ar[d]^{j} & Z \ar[r] \ar[d] & A\ar[r]^{\phi} \ar@{=}[d] &Y[1] \ar[d]^{j[1]} \\
B \ar[r] & X \ar[r]^{i} & A \ar[r]  &B[1] }$$
such that $j\in{{\mcI}^{\perp}}$. The special $\mcI$-preenvelope is defined dually.

An ideal cotorsion pair in $\mcT$ is a pair of ideals  $(\mcI, \mcJ)$ in $\mcT$ for which $\mcI = {^{\perp}}\mcJ$ and $\mcJ=\mcI^{\perp}.$  The ideal cotorsion pair $(\mcI, \mcJ)$ is \emph{complete} if $\mcI$ is a special
precovering ideal and $\mcJ$ is a special preenveloping ideal.

The following result is a direct consequence of Example 3.4.2 and Theorem 5.3.4 in \cite{Breaz}.
\begin{thm} \label{theorem2} Let $\mcI$ be an ideal in $\mcT$. Then the following statements are equivalent.
\begin{enumerate}
\item  $\mcI$ is special precovering in $\mcT$.
\item  $\mcI$ is precovering in $\mcT$.
\item  $(\mcI,\mcI^{\perp})$ is a complete cotorsion pair.
\end{enumerate}
\end{thm}

\section{\bf Proof of Theorem \ref{ThA}}

Throughout this section, $(\mcA; \mcE)$ is a Frobenius category, $\omega$ is the full subcategory consisting of projective objects (also injective objects), and $\mcI$ is an ideal in $(\mcA; \mcE)$ such that $\omega$ $\subseteq$ {\rm Ob}{\rm(}$\mathcal I${\rm)}.

Next we recall the triangulated structure of $\mcA/\omega$ which is induced by a Frobenius category $\mcA$~(see~\cite{Happel}).

Let $\mcA$ be a Frobenius category. For any $A\in \mcA$, we have a triangle
$\xymatrix{A\ar[r]^{\alpha}&I\ar[r]^{\beta}&D\ar@{-->}[r]^{\delta}&,}$
where $I\in \omega$. Define ${T}(A):= D$ as the image of $D$ in $\mcA/\omega$. Then  $[1]:= {T}:\mcA/\omega\rightarrow \mcA/\omega$ is an equivalence as an additive functor.

Let $u:X\rightarrow Y$ be a morphism in $\mcA$. Then there exists the following commutative diagram
$$\xymatrix@C=30pt@R=30pt{X~\ar@{>->}[r]^{m_{x}}\ar[d]^{u}&I\ar@{->>}[r]^{p_{X}}\ar[d]^{i_{u}}&T(X)\ar@{=}[d]\\
Y~\ar@{>->}[r]^{v}&C_{u}\ar@{->>}[r]^{w}&T(X)}$$
as $I$ is injective. We then have the sequence
$\xymatrix@C=0.6cm{\underline{X}\ar[r]^{\underline{u}}&\underline{Y}\ar[r]^{\underline{v}}&
\underline{C_{u}}\ar[r]^{\underline{w}}&\underline{X}[1]}$
in $\mcA/\omega$.
  Sequences obtained in this way are called \emph{standard triangles}. The \emph{distinguished triangles} in $\mcA/\omega$ are defined as the sequences isomorphic to some standard triangles.

\begin{lem}\label{lemma1} $\mathcal I$ is precovering in $\mcA$ if and only if $\mathcal I$/$\omega$ is precovering in $\mcA$/$\omega$.
\end{lem}
\begin{proof} We only prove the ``if" statement, and proof of the ``only if" statement is straightforward. Assume that $\mathcal I/\omega$ is precovering in $\mcA/\omega$. Let ${C}$ be an object in $\mcA$. Then $\underline{C}$ has an $\mathcal I/\omega$-precover $\underline{\alpha} : \underline{X} \to \underline{C}$.
Choose an object ${X_C} \in \mathcal C$ and a morphism ${f_C}$ : ${X_C}$ $\to$ ${C}$ in $\mathcal I$ such that \underline{${X}$$_{C}$} = \underline{${X}$} and \underline{${f}$$_C$} = \underline{$\alpha$}.
Without loss of generality, we may assume that ${f_C}$ is a deflation (if it is not then ${f_C^{'}}$ :=$(\mu, {f_C})$ : ${X'_C}$ := ${T}$ $\oplus$ ${X}$$_{C}$ $\to$ $ C$ is a deflation, where $\mu$ : $ T$ $\to$ $ C$ is a deflation of $ C$ with $T\in{\omega}$, \underline{${X}$$_{C}^{'}$} = \underline{${X}$} and \underline{${f}$$_C^{'}$} = \underline{$\alpha$}).
If ${g}$ : ${X}$${'}$ $\to$ ${C}$ is a morphism in $\mathcal I$, then there exists a morphism \underline{${h}$} : \underline{${X}{'}$} $\to$ \underline{${X}$} in $\mathcal{I/\omega}$ such that \underline{${g}$} = \underline{$\alpha$} \underline{${h}$}. Then ${g-f_C{h}}$ factors through an object $\omega_C$ in $\omega$, i.e. there exists $\lambda$ :${X{'}}$ $\to$ $\omega_C$ and $\phi_C$ : $\omega_C$ $\to$ ${C}$ such that ${g-f_C{h}}$ = $\phi_C{\lambda}$.
Since ${f_C}$ is a deflation, there exists $\sigma$ : $\omega_C$ $\to$ ${X_C}$ such that $\phi_C$ = ${f_C{\sigma}}$. It follows that $g = f_C{(h+\sigma\lambda)}$. Hence ${g}$ factors through ${f_C}$, and this shows that ${f_C}$ is an $\mathcal I$-precover of $C$.
\end{proof}

\begin{lem}\label{lemma2} If $f : X \to\ Y$ is a morphism  in $\mcA$ such that $\underline{f} \in \mathcal I/\omega$, then $f \in \mathcal{I}$.
\end{lem}
\begin{proof} Let $f : X \to Y$ be a morphism in $\mathcal C$ such that $\underline{f} \in \mathcal I/\omega$. Then there exists a morphism $g : X \to Y$ in $\mathcal I(X,Y)$ such that $\underline{f} = \underline{g}$ in $\mathcal I/\omega$. Then $f - g$ factors through an object $\omega_1$ in $\omega$, and therefore there exist $s : X \to \omega_1$ and $t : \omega_1 \to Y$ such that $f - g = t s$. Hence $f \in \mathcal{I}$.
\end{proof}

\begin{lem}\label{lemma3} Let $\mathcal{S}$ be a class of morphisms in $(\mcA; \mcE)$ such that $\omega$ $\subseteq$ {\rm Ob}{\rm(}$\mathcal S${\rm)}. Consider the following homotopy cartesian diagram in $\mcA/\omega$
$$\xymatrix@C=30pt@R=30pt{\underline{Y} \ar[r] \ar@{=}[d] & \underline{C_{1}} \ar[r] \ar[d]^{\underline{j}} & \underline{X_{1}} \ar[r] \ar[d]^{\underline{i}} &\underline{Y}[1] \ar@{=}[d] \\
\underline{Y}\ar[r] & \underline{C} \ar[r]^{\underline{l}} & \underline{X} \ar[r] &\underline{Y}[1],\\ }\eqno{\raisebox{-4ex}{{\rm (3.1)}}}$$
where $i:X_{1}\rightarrow X$, $l:C\rightarrow X$ and $j:C_{1}\rightarrow C$ are morphisms in $\mcA$.
\begin{enumerate}
\item  If $\underline{i}$ belongs to $\mathcal S/\omega$, there exists a pullback diagram
$$\xymatrix{
Y'~\ar@{>->}[r]\ar@{=}[d]&C_{2} \ar[d]\ar@{->>}[r]&X_{1}\ar[d]^{i}\\
Y'~\ar@{>->}[r]&C' \ar@{->>}[r]&X,\\}\eqno{\raisebox{-4ex}{{\rm (3.2)}}}$$
 such that the induced diagram
$$\xymatrix@C=30pt@R=30pt{\underline{Y'} \ar[r] \ar@{=}[d] & \underline{C_{2}} \ar[r] \ar[d] & \underline{X_{1}} \ar[r] \ar[d]^{\underline{i}} &\underline{Y'}[1] \ar@{=}[d] \\
\underline{Y'} \ar[r] & \underline{C'} \ar[r] & \underline{X} \ar[r] &\underline{Y'}[1]\\}$$
is isomorphic to the given bicastesian diagram (3.1) in $\mathcal{C}/\omega$.

\item If $\underline{j}$ belongs to $\mathcal S/\omega$, there exists a pushout diagram
$$\xymatrix{
C_{1}~\ar@{>->}[r]\ar[d]^{j}&X_{3} \ar[d]\ar@{->>}[r]&Y_{1}\ar@{=}[d]\\
C~\ar@{>->}[r]&X' \ar@{->>}[r]&Y_{1},\\}\eqno{\raisebox{-4ex}{{\rm (3.3)}}}$$
such that the induced diagram
$$\xymatrix@C=30pt@R=30pt{\underline{C_{1}} \ar[r] \ar[d]^{\underline{j}} & \underline{X_{3}} \ar[r] \ar[d] & \underline{Y_{1}} \ar[r] \ar@{=}[d] &\underline{C_{1}}[1] \ar[d] \\
\underline{C} \ar[r] & \underline{X'} \ar[r] & \underline{Y_{1}} \ar[r]  &\underline{C_{3}}[1]\\ }$$
is isomorphic to the given bicastesian diagram (3.1) in $\mathcal{C}/\omega$.
\end{enumerate}
\end{lem}

\begin{proof} We only prove (1), and the proof of (2) is similar. It suffices to construct the desired pullback diagram (3.2). Note that there exists a deflation $p:P\rightarrow X$ with $P$ a projective object in $\mcA$. It follows that $(l,p):C\oplus P\rightarrow X$ is also a deflation. If we set $C':=C\oplus P$, there exists a conflation $Y'\rightarrow C'\rightarrow X$ in  $\mcA$. This yields the desired pullback diagram  (3.2).
\end{proof}

The {\em arrow category} $\Arr (\mcA)$ of a category $\mcA$ is the category whose objects \linebreak $a: A_0 \to A_1$ are the morphisms (arrows) of $\mcA,$
and a morphism $f: a \to b$ in $\Arr (\mcA)$ is given by a pair of morphisms $f = (f_0, f_1)$ in $\mcA$ for which the diagram
$$
\xymatrix@C=30pt@R=30pt{A_0 \ar[r]^{f_0} \ar[d]^{a} & B_0 \ar[d]^b \\
A_1 \ar[r]^{f_1} & B_1}$$
commutes. In addition, the additive category $\Arr (\mcA)$ may be equipped with an exact structure for which the collection $\Arr (\mcE)$ of distinguished kernel-cokernel pairs
$$\xi: \xymatrix@C=30pt@R=30pt{b~\ar@{>->}[r]^m & c \ar@{->>}[r]^p & a,}$$ is given by
morphisms of conflations in $(\mcA; \mcE),$
$$
\xymatrix@C=30pt@R=35pt{\Xi_0: B_0~ \ar@<-2.5ex>[d]_{\xi} \ar@<2ex>[d]^b \ar@{>->}[r]^-{m_0} & C_0 \ar@{->>}[r]^{p_0} \ar[d]^c & A_0 \ar[d]^a \\
\Xi_1: B_1~ \ar@{>->}[r]^-{m_1} & C_1 \ar@{->>}[r]^{p_1} & \; A_1.}$$

Recall from \cite{FH} that a conflation of arrows $\xi: \xymatrix@C=30pt@R=30pt{b~\ar@{>->}[r] & c \ar@{->>}[r] & a}$ in the exact structure $(\Arr (\mcA); \Arr (\mcE))$ is said to be $\ME$ ({\em mono-epi}) if there exists a factorization
$$\xymatrix@C=30pt@R=30pt{
\Xi_0 : B_0~ \ar@{>->}[r] \ar@{=}@<2ex>[d] \ar@<-3ex>[dd]_{\xi} & C_0 \ar@{->>}[r] \ar[d]^{c_1} & A_0 \ar[d]^a \\
 \;\; \;\;\;\; B_0~ \ar@{>->}[r] \ar@<2ex>[d]^b & C \ar@{->>}[r] \ar[d]^{c_2} & A_1 \ar@{=}[d]  \\
\Xi_1 : B_1~ \ar@{>->}[r] & C_1 \ar@{->>}[r] & A_1,
}$$ with $c = c_2 c_1.$

\begin{prop}\label{prop-cotorsion-pair} Let $\mathcal{S}$ be a class of morphisms in $(\mcA; \mcE)$ such that $\omega$ $\subseteq$ {\rm Ob}{\rm(}$\mathcal S${\rm)}. Then we have $\mathcal{(S/\omega)^{\perp}} = \mathcal{S^{\perp}/\omega}$ and $^{\perp}{\mathcal{(S/\omega)}} = \mathcal{^{\perp}S/\omega}$.
\end{prop}
\begin{proof} We only prove the first equality, and the proof of the second equality is similar. Let $i : X_1 \to  X$ in $\mathcal S$ and $\underline j : \underline Y \to \underline{Y_1}$ in $\mathcal{(S/\omega)^{\perp}}$. Then there exists a morphism $j' : Y' \to {Y_1}'$ in $\mathcal C$ such that $\underline{j'} = \underline{j}$. For every ME-extension
$$\xymatrix@C=30pt@R=30pt{
 Y'~\ar@{>->}[r]^{f} \ar@{=}[d]  & C_1 \ar@{->>}[r] \ar[d] & X_1 \ar[d]^{i} \\
 Y'~\ar@{>->}[r] \ar[d]^{j'} & C \ar@{->>}[r] \ar[d] &X\ar@{=}[d]  \\
Y'_{1}~\ar@{>->}[r] & C_2 \ar@{->>}[r] & X,}$$
we have the following commutative diagram:
$$\xymatrix@C=30pt@R=30pt{
 \underline{Y'} \ar[r]^{\underline{f}} \ar@{=}[d]  & \underline{C_1} \ar[r] \ar[d] & \underline{X_1} \ar[d]^{\underline{i}} \ar[r]& \underline{Y'}[1] \ar@{=}[d] \\
 \underline{Y'} \ar[r] \ar[d]^{\underline{j'}}  & \underline{C} \ar[r] \ar[d] & \underline{X} \ar@{=}[d] \ar[r]& \underline{Y'}[1] \ar[d] \\
 \underline{Y'_{1}} \ar[r] & \underline{C_2} \ar[r] & \underline{X} \ar[r]& \underline{Y'_{1}}[1].}$$
 Since $\underline{j'} = \underline{j} \in \mathcal{(S/\omega)^\perp}$ and $\underline{i} \in \mathcal{S/\omega}$, there exists a morphism $\underline{g} : \underline{C_1} \to \underline{{Y_1}'}$ such that $\underline{j'} = \underline{g} \underline{f}$. Then $j'-gf$ factors through an object in $\omega$, and therefore there exists a factorization
 \begin{center}
 $j' - g f = t s$ : $Y' \stackrel{s}{\longrightarrow} \omega_1 \stackrel{t}{\longrightarrow} {Y_1}'$,
 \end{center}
 with $\omega_1 \in \omega$. Thus there exists a morphism $k : C_1 \to \omega_1$ such that $s = k f$, and hence $j' = g f + t s = (g + t k) f$ which implies $j' \in \mathcal{S^\perp}$. Since $\underline{j'} = \underline{j}$ and $\underline{j} \in \mathcal{(S/\omega)}^\perp$, we have $\underline{j} \in \mathcal{{S^\perp}/\omega}$, which shows $\mathcal{(S/\omega)}^\perp \subseteq \mathcal{S^\perp}/\omega$.

For the reverse containment $\mathcal{S^\perp}/\omega\subseteq \mathcal{(S/\omega)}^\perp $, we choose $j_1 : Y_2 \to Y_3$ in $\mathcal S^\perp$ and $\underline{i_1} : \underline{X_3} \to \underline{X_2}$ in $\mathcal S/\omega$. We consider the following commutative diagram
$$\xymatrix@C=30pt@R=30pt{
 \underline{Y_{2}} \ar[r] \ar@{=}[d]  & \underline{C_3} \ar[r] \ar[d] & \underline{X_3} \ar[d]^{\underline{i_{1}}} \ar[r]& \underline{Y_{2}}[1] \ar@{=}[d] \\
 \underline{Y_{2}} \ar[r] \ar[d]^{\underline{j_{1}}}  & \underline{C_{2}} \ar[r] \ar[d] & \underline{X_{2}}
 \ar@{=}[d] \ar[r]& \underline{Y_{2}}[1] \ar[d] \\
 \underline{Y_{3}} \ar[r] & \underline{C_4} \ar[r]^{\underline{l}} & \underline{X_{2}} \ar[r]&
 \underline{Y_{3}}[1].}$$
 Thanks to Lemma \ref{lemma3}, we have the following commutative diagram
 $$\xymatrix@C=30pt@R=30pt{
 Y_{2}~\ar@{>->}[r] \ar@{=}[d]  & C'_3\ar@{->>}[r] \ar[d] & X'_3 \ar[d]^{i'_{1}} \\
 Y_{2}~\ar@{>->}[r] \ar[d] & C'_{2} \ar@{->>}[r] \ar[d] &X'_{2}\ar@{=}[d]  \\
Y_{3}~\ar@{>->}[r] & C'_4 \ar@{->>}[r]^{l'} & X'_{2},}$$
with $\underline{{i'_{1}}}$ = $\underline{i_1}$ and $\underline{l'} = \underline{l}$. Since $j_1 \in \mathcal{S^\perp}$, there exists a morphism $s' : {X'_3}\to {C'_4}$ such that ${i'_1} = l's'$. Then we have $\underline{i_1} = {\underline{l}}  \underline{s'}$, and thus we have $\underline{j_1} \in \mathcal{(S/\omega)^\perp}$. This completes the proof.
\end{proof}

\begin{prop}\label{thm1} The following statements are true:
\begin{enumerate}
\item $\mcI$ is special precovering in $\mcA$ if and only if $\mathcal I/\omega$ is special covering in $\mcA/\omega$.
\item $\mcI$ is special preenveloping in $\mcA$ if and only if $\mathcal I/\omega$ is special preenveloping in $\mcA/\omega$.
\end{enumerate}
\end{prop}
\begin{proof} We only prove (1), and the proof of (2) is similar. The ``only if" statement follows from Lemma \ref{lemma1} and \cite[Example 3.4.2]{Breaz}. For the ``if" statement, we assume that
$\mathcal I/\omega$ is special covering in $\mcA/\omega$. Let $C$ be an object in $\mcA$. Then there exists an $\mathcal I/\omega$-precover $\underline{i} : \underline{X} \to \underline{C}$ and a homotopy cartesian diagram
$$\xymatrix@C=30pt@R=30pt{\underline{Y} \ar[r]^{\underline{h}} \ar[d]^{\underline{j}} & \underline{Z} \ar[r] \ar[d]^{\underline{t}} & \underline{C} \ar[r] \ar@{=}[d] &\underline{Y}[1] \ar[d] \\
\underline{B}\ar[r]^{\underline{f}} & \underline{X} \ar[r]^{\underline{i}} & \underline{C} \ar[r] &\underline{B}[1],\\ }\eqno{\raisebox{-4ex}{{\rm (3.4)}}}$$
with $\underline{j} \in \mathcal {(I/\omega)^\perp}$. Thanks to Lemma \ref{lemma1}, there exists a conflation $ B' \stackrel{f'}\mto X_C \stackrel{i'}\eto C$ satisfying that $i'$ is an $\mathcal I$-precover of $C$, $\underline{X_C} = \underline{X}$, $\underline{i'} = \underline{i}$ and $\underline{B'} \stackrel{\underline{f'}}{\longrightarrow} \underline{X_C} \stackrel{\underline{i'}}{\longrightarrow} \underline{C} \stackrel{}\longrightarrow \underline{B'[1]}$ is isomorphic to $\underline{B} \stackrel{\underline{f}}{\longrightarrow} \underline{X} \stackrel{\underline{i}}{\longrightarrow} \underline{C} \stackrel{}\longrightarrow \underline{B[1]}$.
Since $\underline{t} \in \mathcal{C}(Z,X)/\omega(Z,X)$ and $\underline{X} = \underline{X_C}$, there exists a deflation $t' : Z' \to X_C$ such that $\underline{t} = \underline{t'}$.

Let $m := i' t{'}$. Then $m:Z'\rightarrow C$ is a deflation. Thus there exists a conflation $Y'\s{l}\rightarrow Z'\s{m}\rightarrow C$ in $\mcA$, whence we obtain the following commutative diagram
$$\xymatrix@C=30pt@R=30pt{
Y'\ar@{>->}[r]^{l}\ar[d]^{k}&Z' \ar[d]^{t'}\ar@{->>}[r]^{m}&C\ar@{=}[d]\\
B'\ar@{>->}[r]^{f'}&X_{C} \ar@{->>}[r]^{i'}&C\\}$$
 in $\mcA$, which induces the following homotopy
cartesian diagram
$$\xymatrix@C=30pt@R=30pt{
\underline{Y'}\ar[r]^{\underline{l}}\ar[d]^{\underline{k}}&\underline{Z'} \ar[d]^{\underline{t'}}\ar[r]^{m}&\underline{C}\ar@{=}[d]\ar[r]&\underline{Y'}[1]\ar[d]\\
\underline{B'}\ar[r]^{\underline{f'}}&\underline{X_{C}} \ar[r]^{\underline{i'}}&\underline{C}\ar[r]&\underline{B'}[1]\\}\eqno{\raisebox{-4ex}{{\rm (3.5)}}}$$
 in $\mcA/\omega$. It is clear that two homotopy
cartesian diagrams  (3.4) and  (3.5) are isomorphic. This implies that $\underline{k}$ belongs to $ \in \mathcal {(I/\omega)^\perp}$. Thanks to Lemma \ref{lemma2}, we have ${k}\in{\mathcal{I}^{\perp}}$. So $i':X_{C}\rightarrow C$ is a special $\mathcal{I}$-precover of $C$, as desired.
\end{proof}
We are now in a position to prove Theorem \ref{ThA} in the introduction.

{\bf Proof of Theorem \ref{ThA}.} We only need to prove (1), and the proof of (2) is similar. Assume that $\mcI$ is an ideal in $(\mcA; \mcE)$ such that $\omega$ $\subseteq$ {\rm Ob}($\mathcal I$). $(a)\Leftrightarrow(b)$ is trivial.

$(a)\Leftrightarrow(c)$ holds by \cite[Theorem 1]{FGHT}.

$(c)\Leftrightarrow(d)$ follows from Proposition \ref{thm1}.

$(b)\Leftrightarrow(d)$ follows from Lemma \ref{lemma1} and \cite[Example 3.4.2 and Theorem 5.3.4]{Breaz}.  \hfill$\Box$
\vspace{2mm}







\section{\bf Applications of Theorem \ref{ThA}}\label{application:Enochs-conjecture}
Assume that $(\mcA; \mcE)$ is an exact category with exact coproducts. A celebrated Bongartz-Eklof-Trlifaj Lemma \cite[Theorem 13]{FHHZ} says that the ideal $a^{\perp}$ is special preenveloping for any morphism $a$ in $\mcA$, which can be seen as the ideal version of the famous Eklof-Trlifaj Lemma \cite[Theorem 2]{ET}. In combination this with Corollary \ref{cor3}, we have the following version of Bongartz-Eklof-Trlifaj Lemma in $\mcA/\omega$.

\begin{cor} \label{cor:BKT-lemma} Let $(\mcA; \mcE)$ be a Frobenius category with exact coproducts and $\underline{a}:\underline{B}\to \underline{C}$ a morphism in $\mcA/\omega$. Then ${\underline{a}}^{\perp}$ is special preenveloping in $\mcA/\omega$. Moreover, $(^{\perp}(\underline{\mathcal{M}}^{\perp}),\underline{\mathcal{M}}^{\perp})$ is a complete cotorsion pair in $\mcA/\omega$ provided that $\mathcal{M}$ is a set of morphisms in $\mathcal{A}$.
\end{cor}
\begin{proof} Thanks to Proposition \ref{prop-cotorsion-pair}, we have
    $\underline{a}^{\perp}=\underline{a^{\perp}}$, $\underline{\mathcal{M}}^{\perp}=\underline{\mathcal{M}^{\perp}}$ and $\underline{^{\perp}(\mathcal{M}^{\perp})}={^{\perp}({\underline{\mathcal{M}^{\perp}}})}$ in $\mcA/\omega$. Thus the result holds by \cite[Theorem 13 and Corollary 15]{FHHZ} and Corollary \ref{cor3}.
\end{proof}

Let $(\mcA; \mcE)$ be an exact category. Recall from \cite{A11} that a pair of subcategories $(\mathcal{X},\mathcal{Y})$ of $\mathcal{A}$ is said to be a \emph{cotorsion pair} if
$\mathcal{X}={^{\perp}\mathcal{Y}}:=\{X\in{\mathcal{A}} \ |\ \Ext_{\mathcal{A}}^{1}(X,Y)=0 \ \textrm{for each} \ Y\in{\mathcal{Y}} \}$ and
$\mathcal{Y}={\mathcal{X}^{\perp}}:=\{Y\in{\mathcal{A}} \ |\ \Ext_{\mathcal{A}}^{1}(X,Y)=0 \ \textrm{for each} \ X\in{\mathcal{X}} \}$. Morover, a cotorsion pair $(\mathcal{X}, \mathcal{Y})$ is called \emph{complete} \cite{A11} if for each $M\in{\mathcal{A}}$, there exist conflations $Y_{M}\mto X_{M}\s{f_{M}}\eto M \ \textrm{and}\ M\s{g^{M}}\mto Y^{M}\eto X^{M}$ in $\mathcal{A}$
such that $X_{M},X^{M}\in{\mathcal{X}}$ and $Y_{M},Y^{M}\in{\mathcal{Y}}$. In this case, $f_{M}:X_{M}\rightarrow M$ is called a \emph{special $\mathcal{X}$-precover} of $M$ and $g^{M}:M\to
Y^{M}$ is called a \emph{special $\mathcal{Y}$-preenvelope} of $M$.

\vspace{2mm}
We are now in a position to prove Corollary \ref{corollary2} in the introduction.


 {\bf Proof of Corollary \ref{corollary2}.}
$(1)\Rightarrow(2)$ holds by \cite[Proposition 3.3]{A11}.

$(2)\Rightarrow(3)$. Let $M$ be an object in ${\rm smd}(\mcX\oplus\mcY)$. It follows that $\underline{M}$ belongs to  ${\rm smd}(\underline{\mcX}\oplus\underline{\mcY})$, and therefore there is a triangle
$$\xymatrix{\underline{M}\ar[r]&\underline{Y}\ar[r]&\underline{X}\ar[r]&\underline{M}[1]}\eqno{{{\rm (4.1)}}}$$  in $\mcA/\omega$ with $Y\in{\mcY}$ and $X\in{\mcX}$ by (2). Thus there exists a conflation $ \xymatrix@C=30pt@R=30pt{M'~\ar@{>->}[r] & Y'\ar@{->>}[r] &X'}$ in $(\mcA; \mcE)$ such that the induced triangle   $\xymatrix{\underline{M'}\ar[r]&\underline{Y'}\ar[r]&\underline{X'}\ar[r]&\underline{M'}[1]}$ is isomorphic to the triangle (4.1) in $\mcA/\omega$. Consequently, there exist projective objects $P$ and $P'$ such that $Y\oplus P\cong{Y'\oplus P'}$. Since $Y\in{\mcY}$ and $\mcY$ is closed under direct summands, so is $Y'$. Similarly, we can show  $X'\in{\mcX}$. Note that $\underline{M}\cong{\underline{M'}}$ in $\mcA/\omega$. Then there exist projective objects $Q$ and $Q'$ such that $M\oplus Q\cong{M'\oplus Q'}$. Consider the following pushout diagram in $\mcA$:
$$\xymatrix@C=30pt@R=30pt{M \ar[r]\ar@{=}[d] & M'\oplus Q'\ar[r] \ar[d] & Q \ar[d] \\
M \ar[r] & Y'\oplus{Q'} \ar[r] \ar[d] & Z\ar[d]\\
&X'\ar@{=}[r] &X',
}$$
where all rows and columns are conflations. Since $Q\in{\omega}$ and $X'\in{\mcX}$, so is $Z$. So the second column in the above diagram yields the desired special $\mcY$-preenvelope of $M$.

$(3)\Rightarrow(1)$. Note that $\mcI=\mcI(\mcX)$ and $\mcJ=\mcI(\mcY)$. Since $(\mcI,\mcJ)$ is a complete ideal cotorsion pair by hypothesis, it follows that $(\mcX,\mcY)$ is a cotorsion pair  such that $\mcX$ is precovering and $\mcY$ is preenveloping in $\mcA$.  Then  $\mcX/\omega$ is precovering and $\mcY/\omega$ is preenveloping in $\mcA/\omega$ by \cite[Lemma 4.1.1]{Beligiannis}.

Consequently, to prove (1), it suffices to show that each object in $\mcA$ has a special $\mcX$-precover. Let $M$ be an object in $\mcA$.
Then $M$ has an $\mcX$-precover $f:X\to M$.
Without loss of generality, we may assume that ${f}$ is a deflation and $K:=\ker(f)$ (if it is not then ${f^{'}}$ :=$(\lambda, {f})$ : ${X'}$ := ${T}$ $\oplus$ ${X}$ $\to$ $ M$ is a deflation, where $\lambda$ : $ T$ $\to$ $ M$ is a deflation of $M$ with $T\in{\omega}$
).
It follows that ${f}:{X}\to{M}$ is an $\mcI$-precover in $\mcA$,
whence we have a triangle in $\mcA/\omega$
$$\xymatrix{\underline{K}\ar[r]^{\underline{k}}&\underline{X}\ar[r]^{\underline{f}}&\underline{M}\ar[r]^{\underline{g}[1]}&\underline{K}[1]}$$
Similar to the proof of Lemma \ref{lemma1}, $\underline{f}$ is an $\mcI/\omega$-precover of $\underline{M}$. By \cite[Example 3.4.2]{Breaz} and Corollary \ref{corollary2}, $\underline{g}[1]\in\mcJ/\omega[1]$.

On the other hand, $\underline{M}[-1]$ has a ${\mcY}$-preenvelope $\alpha:\underline{M}[-1]\to \underline{Y}$ in $\mcA/\omega$ with $Y\in{\mcY}$. Thus, $\alpha:\underline{M}[-1]\to \underline{Y}$ is an $\mcJ/\omega$-preenvelope, and hence we have a triangle in $\mcA/\omega$
$$\xymatrix{\underline{Y}\ar[r]^{\underline{s}}&\underline{L}\ar[r]^{\underline{\beta}}&\underline{M}\ar[r]^{\underline{\alpha}[1]}&Y[1]}.$$
By \cite[Example 3.4.2]{Breaz} and Corollary \ref{corollary2} again, $\underline{\beta}\in\mcI/\omega$ and $\underline{\alpha}[1]$ is an $\mcJ/\omega[1]$-preenvelope of $\underline{M}$. Note that $\Hom_{\mcA/\omega}(\underline{X},\underline{Y}[1])=0$ and $\underline{g}[1]\in\mcJ/\omega[1]$,  then we obtain the following commutative diagram in $\mcA/\omega$:
$$\xymatrix@C=30pt@R=30pt{
 \underline{K} \ar[r]^{\underline{k}}
 \ar@{-->}[d]^{\underline{\varphi}}  &  \underline{X}\ar[r]^{\underline{f}} \ar@{-->}[d] & \underline{M} \ar@{=}[d] \ar[r]^{\underline{g}[1]} &\underline{K}[1]\ar@{-->}[d]^{\underline{\varphi}[1]}\\
 \underline{Y} \ar[r]^{\underline{s}}
 \ar@{-->}[d]^{\underline{\psi}}  &  \underline{L}\ar[r]^{\underline{\beta}} \ar@{-->}[d] & \underline{M} \ar@{=}[d] \ar[r]^{\underline{\alpha}[1]} &\underline{Y}[1]\ar@{-->}[d]^{\underline{\psi}[1]}\\
\underline{K} \ar[r]^{\underline{k}}  &  \underline{X}\ar[r]^{\underline{f}} & \underline{M} \ar[r]^{\underline{g}[1]} &\underline{K}[1].}$$
It follows that $\underline{\psi}[1]\underline{\alpha}[1]=\underline{g}[1]$ and $\underline{\varphi}[1]\underline{g}[1]=\underline{\alpha}[1]$. Thus we have $(1-\underline{\psi}[1]\underline{\varphi}[1])\underline{g}[1]=0$, whence there exists $\underline{\mu}[1]:\underline{X}[1]\to \underline{Y}[1]$ such that $\underline{\mu}[1]\underline{k}[1]=1-\underline{\psi}[1]\underline{\varphi}[1]$. This implies that $(\underline{k}[1],\underline{\varphi}[1])^{T}:\underline{K}[1]\rightarrow \underline{X}[1]\oplus \underline{Y}[1]$ is a split monomorphism in $\mcA/\omega$, and therefore $\underline{K}$ is isomorphic to a direct summand of $\underline{X}\oplus \underline{Y}$ in $\mcA/\omega$. Thus there exist an object $K'$ in $\mcA$, projective objects $P$ and $Q$ such that $K\oplus K'\oplus P\cong X\oplus Y\oplus Q$, whence $K$ belongs to ${\rm smd}(\mcX\oplus\mcY)$. By hypothesis, $K$ has a special $\mcY$-preenvelope $\upsilon:K\to V$. Consider the following pushout diagram in $\mcA$:
$$\xymatrix@C=30pt@R=30pt{K\ar[r]\ar[d]^{\upsilon} & X\ar[r]^{f} \ar[d] & M \ar@{=}[d] \\
V \ar[r]\ar[d] & Z \ar@{->>}[r] \ar[d] & M \\
U\ar@{=}[r] &U,
}$$
where all rows and columns are conflations. Since $U$ and $X$ are in $\mcX$, so is $Z$. So the second row in the above diagram yields the desired special $\mcX$-precover of $M$.

For the second part, we choose an object $\underline{M}$  in ${\rm smd}(\underline{\mcX}\oplus\underline{\mcY})$. Cleary, $\underline{\mcY}$ is preenveloping in $\mcA/\omega$. Since $\mcA/\omega$ is a Krull-Schmidt triangulated category by hypothesis, there is a triangle
$\xymatrix{\underline{M}\ar[r]^{\underline{a}}&\underline{C}\ar[r]&\underline{F}\ar[r]&\underline{M}[1]}$  in $\mcA/\omega$ with $C\in{\mcY}$ such that $\underline{a}$ is a minimal left $\underline{\mcY}$-approximation, that is, $\underline{a}$ is a $\underline{\mcY}$-preenvelope and $\underline{a}$ does not
have a direct summand of the form $0\to \underline{T}$
$(\underline{T}\in{\mcA/\omega}, \underline{T}\neq 0)$ as a complex (see \cite{IY}). Since $\underline{\mcX}=\{\underline{M}\in{\mcA/\omega}~|~ \Hom_{\mcA/\omega}(\underline{M},\underline{Y}[1])=0~\textrm{for}~\textrm{any}~\underline{Y}\in
{\underline{\mcY}}\}$ by \cite[Proposition 3.3]{A11}, it follows that $\underline{F}$ belongs to $\underline{\mcX}$ and the proof is dual to Proposition 2.3(1) in \cite{IY}. This shows $F\in{\mcX}$, as desired.\hfill$\Box$

\vspace{4mm}
\renewcommand\refname

\vspace{4mm}
\small

\noindent\textbf{Dandan Sun}\\
School of Mathematics Sciences, Zhejiang University of Technology, Hangzhou 310023, China\\
E-mail: 13515251658@163.com\\[1mm]
\textbf{Zhongsheng Tan}\\
School of Mathematics Sciences, Zhejiang University of Technology, Hangzhou 310023, China\\
E-mails: 3186856595@qq.com\\[1mm]
\textbf{Qikai Wang}\\
School of Mathematics Sciences, Zhejiang University of Technology, Hangzhou 310023, China\\
E-mails: qkwang@zjut.edu.cn\\[1mm]
\textbf{Haiyan Zhu}\\
School of Mathematics Sciences, Zhejiang University of Technology, Hangzhou 310023, China\\
E-mails: hyzhu@zjut.edu.cn\\[1mm]
\end{document}